\documentclass[11pt,reqno]{amsart}
\usepackage{amsfonts}
\usepackage{fancyhdr}
\usepackage{times}
\usepackage[T1]{fontenc}
\usepackage{mathrsfs}
\usepackage{latexsym}
\usepackage[dvips]{graphics}
\usepackage{epsfig}
\usepackage{hyperref, amsmath, amsthm, amsfonts, amscd, flafter,epsf}
\usepackage{amsmath,amsfonts,amsthm,amssymb,amscd}
\input amssym.def
\input amssym.tex
\usepackage{color}
\usepackage[parfill]{parskip}
 \usepackage{mathrsfs}
\usepackage{geometry}              
\usepackage{epstopdf}

\usepackage{verbatim}

\newtheorem{thm}{Theorem}

\newtheorem{con}[thm]{Conjecture}
\newtheorem{lem}[thm]{Lemma}
\newtheorem*{thm*}{Theorem}
\newtheorem*{con*}{Conjecture}
\newtheorem*{lem*}{Lemma}

\newtheorem{defn}{Definition}[section]
\newtheorem{rem*}{Remark}

\begin{document}
 \baselineskip=17pt
\hbox{}
\medskip

\title
[solution to index conjecture] {Solution to the index conjecture in zero-sum theory}

\author{Fan Ge}

\email{fge@wm.edu}

\address{Department of Mathematics, William \& Mary, Williamsburg, VA, United States}


\maketitle

\begin{abstract}
	A problem in zero-sum theory is to determine all pairs $(k,n)$ for which every minimal zero-sum sequence of length $k$ modulo $n$ has index $1$. While all other cases have been solved more than a decade ago, the case when $k$ equals $4$ and $n$ is coprime to $6$ remains open. Precisely, The Index Conjecture in this subject states that if $n$ is coprime to $6$ then every minimal zero-sum sequence of length $4$ modulo $n$ has index $1$. In this paper, we prove an equivalent version of this conjecture for all $n>N$ for some absolute constant $N$. 
\end{abstract}


\medskip

\section{Introduction}

Throughout this paper let $G=\mathbb Z/n$ be an additive cyclic group of order $n$.  By a \emph{sequence $S$ of length $k$ over $G$}
we mean a sequence with $k$ elements, each of which is in $G$. We
write $(a_1)\cdots(a_k)$ for such a sequence, rather than $a_1$,
\dots, $a_k$. A sequence $S$ is said to be a \emph{zero-sum}
sequence if $\sum_i a_i=0$. It is a \emph{minimal} zero-sum sequence
if it is a zero-sum sequence but no proper nontrivial subsequence of
it is. 

For each element $a$ in $G$, one can identify $a$ with an integer in the interval $[0,n-1]$. Further, for each integer $x$ we denote by $(x)_n$ the least non-negative integer in the residue class of $x$ modulo $n$. We also require the following definition. 

\begin{defn}\textnormal{For a sequence over $G$
		$$S=(a_1)\cdots(a_k), \qquad \textnormal{where} \ 1\le a_1,...,a_k\le n\,,$$
		define the \textit{index} of $S$, written as $\text{ind}(S)$, to be smallest number in the set
		$$\left\{\frac{1}{n}\cdot \sum_{j=1}^{k}(ga_j)_n: g\in (\mathbb{Z}/n)^*\right\}$$ }
where $(\mathbb{Z}/n)^*$ is the multiplicative group in $ \mathbb{Z}/n $.
\end{defn}

The index of a sequence is a crucial invariant in the theory of zero-sum sequences. The study of indices (in different terminology) dates back to 1980s (see~\cite{KL, Ger90}). It was first defined and given its name by Chapman et al.~\cite{CFS} in 1999.
It has been used as an important tool in the
study of zero-sum sequences and related topics (see, for example, \cite{Ger09, Gry13}).
The study of indices also finds its connection with integer partitions~\cite{P}, factorization theory~\cite{Gao-Ger,Ger09}, Dedekind sums~\cite{LPYZ, Ge1}, and Heegaard
Floer Homology~\cite{JRW, Ge1}, as well as discrepancy estimates which is pointed out to me by Granville~\cite{Gran}.   A fundamental question is to determine the pairs
$(k, n)$ for which every minimal zero-sum sequence $S$ of length $k$
over $G$ has index $1$. We shall call such pairs $(k,n)$ \textit{good} pairs. When $k\le 3$ every pair $(k,n)$ is trivially good. When $5\le k\le \frac{n}{2}+1$ it is known that every pair is bad; This result has been in literature for at least 10 years, but the author was not able to identify who first proved it. When $k>\frac{n}{2}+1$ we know each pair is good; This result was proved independently by Savchev and Chen~\cite{S-C}, and by Yuan~\cite{Y}. (See also Gao~\cite{Gao}). The case $k=4$ is more subtle. When $k=4$ but $\gcd(n, 6)\ne 1$, Ponomarenko~\cite{P} showed that there are sequences with index greater than $1$. He also suggested the following conjecture, and verified it for $n\le 1000$.
\begin{con*}
\textnormal{(The Index Conjecture, Version 1)} Let $\gcd(n, 6) = 1$. Then
	every minimal zero-sum sequence $S$ over $G$ of length $4$ has
	\textnormal{ind}($S$) = 1.
\end{con*}

The study of the index conjecture has been active~\cite{LPYZ, LP, XS, SXL, SX, X, Ge1,Ge2, ZQ, GV}, and we have a number of partial results. In particular, a consequence of a result in Shen et al.~\cite{SXL}  is (see~\cite{Ge1})
\begin{thm*}
	Suppose that $n$ is the smallest integer for which
	Conjecture~\ref{conj} fails. Let $S=(a_1)(a_2)(a_3)(a_4)$ be a
	minimal zero-sum sequence over $G=\mathbb{Z}/n$ with
	\textnormal{ind}$(S)=2$. Then we have $\gcd(n, a_i)=1$ for all $i$.
\end{thm*}
Therefore, to prove the index conjecture, it suffices to prove the following variant. 
\begin{con}\label{conj}
	\textnormal{(The Index Conjecture, Version 2)}	Let $\gcd(n, 6) = 1$. Let $S=(a_1)(a_2)(a_3)(a_4)$ be a minimal zero-sum sequence over $G$. Suppose $\gcd(n, a_i)=1$ for all $i$. Then
	\textnormal{ind}($S$) = 1.
\end{con}

Recently, Zeng and Qi~\cite{ZQ} proved that the index conjecture is true if $\gcd(n,30)=1$. In a more recent preprint~\cite{GV} Grynkiewicz and Vishne gave an alternative proof of some previous results. They also wrote:
\begin{quote}
\begin{center}``As will be seen in the proof, the majority of the complications that arise in our arguments happen only
when $5\mid n$ (more than half our arguments are devoted solely to this case), giving further indication
that the case when $5\mid n$ is fundamentally harder than the case $gcd(n, 30) = 1$.''\end{center}
\end{quote}

In this paper, we settle the index conjecture completely, except for a finite number of $n$'s.
\begin{thm}\label{thm index}
	Conjecture~\ref{conj} is true for all $n>N= 10^{20}$.
\end{thm}

Thus, to prove the index conjecture, it  remains to verify Conjecture~\ref{conj} for $n\le N.$ I have not attempted to obtain the optimal bound $N$.

\section{Notation}
Our proof for Theorem~\ref{thm index} is Fourier analytic (and of course, number theoretic). It seems that this is the first time Fourier analytic method is introduced in the study of indices of zero-sum sequences. 

Before we present the proof, it is necessary to introduce some notation.
Throughout, $n$ is coprime to $6$. We use  $(k,m)$ to denote the greatest common divisor of two integers $k$ and $m$, unless the context clearly suggests it is an ordered pair. As usual, $e(x)=e^{2\pi i x}$, $e_n(x)=e^{2\pi i x/n}$, and $\phi(\cdot)$ is the Euler totient function. The letter $g$ will denote a generic member in $\left(\mathbb Z/n\right)^*$, and the Ramanujan sum $c_n(m)$ is defined as $\sum_g e_n(mg)$, where the sum runs over all $g$.

We let $I$ be the interval $[0,1/2]$, and $\chi$ be the period-$1$ characteristic function for $I$. That is,
\begin{equation*} \chi(t) = \left\{
\begin{array}{ccc} & 1,  & \textrm{ if }\ 0<\{t\}<1/2,\\[1ex]  & 1/2,  & \textrm{
	if } t\in \frac{1}{2}\mathbb Z,\\[1ex] & 0,  &\textrm{ if }\ \{t\}>1/2.
\end{array}\right.
\end{equation*}
The Fourier coefficients of $\chi$ are defined as
$$\hat\chi(k)=\int_{-1/2}^{1/2} \chi(t)e(-tk)dt.$$
Hence, by a straightforward calculation, we have
\begin{equation*} \hat\chi(k) = \left\{
\begin{array}{ccc} & 1/2,  & \textrm{ if }\ k=0,\\[1ex]  & 0,  & \textrm{
	if } k\ne 0, 2\mid k,\\[1ex] & 1/(i\pi k),  &\textrm{ if }\ 2\nmid k.
\end{array}\right.
\end{equation*}

We will use a smoothed version of $\chi$, denoted by $f$, so that the Fourier coefficients of $f$ are of finite support. For this we follow the idea in~\cite{Vaa} (See also~\cite{D-T, M}). Precisely, we let
$$H(z)=\left(\frac{\sin \pi z}{\pi }\right)^2\left(\sum_{m\in\mathbb Z}\frac{\text{sgn}(m)}{(z-m)^2}+\frac{2}{z}\right)$$
and $J(z)=\frac{1}{2}H'(z).$
We then set $J_{H+1}(x)=(H+1)J((H+1)x)$ where $H$ is a parameter to be chosen later, and define
 $$f(x)=(\chi * j_H)(x)=\int_{-1/2}^{1/2} \chi(x-t)j_H(t)dt,$$
where $j_H(x)=\sum_{m\in\mathbb Z}J_{H+1}(x+m)$.
Therefore, $$f(x)=\sum_h \hat{f}(h)e(hx)=\sum_h \hat{\chi}(h)\hat{j_H}(h)e(hx),$$
where $\hat{j_H}(h)=\hat{J}_{H+1}(h)$, the latter being the Fourier transform of $J_{H+1}$ on $\mathbb R$.
The importance of this choice of $f$ is that the Fourier transform of $J$ has compact support. More precisely, by Lemma 1.23 in~\cite{D-T}, $\hat{J}(t)$ is a real function supported on $[-1,1]$, $| \hat{J}(t)|\le 1$, and $\hat{J}(0)=1$. Hence we conclude that
$$f(x)=\sum_{|h|\le H} \hat{f}(h)e(hx),$$ where 
$\hat{f}(0)=1/2$, $\hat{f}(k)=0$ for even non-zero $k$, and $|\hat{f}(k)|\le 1/(\pi|k|)$ for odd $k$.

We will need an estimate on how well $f$ approximates $\chi$. For this we have (cf.~\cite{D-T}, page 20-21)
$$|\chi(x)-f(x)|\le \frac{1}{H+1}\sum_{|h|\le H}\hat{K}_{H+1}(h)C_he(hx),$$
where $|C_h|\le 1$, and $\hat{K}_{H+1}$ is the Fourier transform of $K_{H+1}(x)=(H+1)K((H+1)x)$ with
$$K(z)=\left(\frac{\sin \pi z}{\pi z}\right)^2$$ being the Fej\'er's kernel. 
We will only use the fact that $|\hat{K}_{H+1}(h)|\le 1.$

Theorem~\ref{thm index} will be a consequence of the following result.
\begin{thm}\label{thm}
	Let $1,a,b$ be integers coprime to $n$. Assume $(n,6)=1$. Suppose that none of the pairs $(1,a),(1,b),(a,b)$ satisfies any relations $x\pm y\equiv 0\pmod n$ for a pair $(x,y)$, and that at least one of the three pairs does not satisfy any relations $x\pm 3y\equiv 0 \pmod n$ or $3x\pm y\equiv 0\pmod n$. 
	Then
	$$\#\{g:g\in \left(\mathbb{Z}/n\right)^*, (g)_n/n\in I,(ag)_n/n\in I,(bg)_n/n\in I\}\ge c\cdot \phi(n)$$ 
	for some absolute constant $c>0$, as long as $n>N$ for some absolute constant $N$. 
\end{thm}


\section{Lemmas}

\begin{lem}\label{lem 3 fall in I}
	Let $(1)(a)(b)(c)$ be a minimal zero-sum sequence modulo $n$ whose index is $2$, and assume $(n,6)=1.$ Then for any $g\in (\mathbb Z/n)^*$ exactly two of
	$(g)_n,(ga)_n,(gb)_n,(gc)_n$ lie in the interval $(0,n/2)$.
\end{lem}
\proof See Remark 2.1 of \cite{LP}.

Lemma~\ref{lem 3 fall in I} suggests a geometric way to think of the index conjecture. Namely, we may view $(1,a,b)$ as a point in $\mathbb R^3$ and consider set of points $$\left\{\left(\frac{(g)_n}{n},\frac{(ga)_n}{n},\frac{(gb)_n}{n}\right):g\in \left(\mathbb{Z}/n\right)^*\right\}.$$
If this set intersects with the cube $I\times I\times I$, then by Lemma~\ref{lem 3 fall in I} the sequence has index $1$. Such a problem is naturally related to discrepancy estimates. A related problem was studied by Granville et al. in~\cite{GSZ}, in which the set of points comes from a 'nondegenerate' curve over a prime field. A key tool in~\cite{GSZ} is Bombieri's deep estimates on exponential sums along curves~\cite{B}. The main difference between the situation in~\cite{GSZ} and ours here is that our set of points is not a curve over finite fields. Indeed, our modulus $n$ is composite in general, instead of a prime number or a prime power, and therefore it seems that Bombieri's result cannot be applied directly. Moreover, it has been known that for the index conjecture, the prime (power) case is significantly simpler than the general case. 

The next four lemmas are devoted to providing, in our case, analogue results of some consequences of Bombieri's square root cancellation estimates.
\begin{lem}\label{lem gcd}
	Let $G$ and $A$ be fixed integers. For any fixed integer $x$, there is at most one integer $y\in[-G,G]$ such that $(xA+y,n)>\sqrt{2Gn}$. Furthermore, if $(A,n)=1$, then for any fixed integer $y$, there is at most one integer $x\in[-G,G]$ such that $(xA+y,n)>\sqrt{2Gn}$.
\end{lem}

\proof For the first statement it suffices to prove that for any two distinct integers $y,z\in[-G,G]$ we have
$$\min\left((xA+y,n),(xA+z,n)\right)\le\sqrt{2Gn}.$$
Note that $(xA+y,xA+z)=(xA+y,y-z)\le |y-z|\le 2G$, thus
$$(xA+y,n)\cdot(xA+z,n)\le 2Gn$$
and the result follows.

If in addition $(A,n)=1$, then $(xA+y,n)=(x+yA^{-1},n)$ where $AA^{-1}\equiv 1 \pmod n.$ Then use the first statement. \qed

\begin{lem}\label{lem *sum}
Let $(A,n)=(n,6)=1$. Let $H$ be a positive integer, and $k$ be an odd integer. Let $k^*$ (depending on $A,n,H$) be the only possible integer in $[-H^2,H^2]$ such that $(Ak+k^*,n)>\sqrt{2H^2n}$.  Let $S$ (also depending on $A,n,H$) denote the set $\{(k,k^*):\ |k|\le H, \ 2\nmid k, \ |k^*|\le H,  \ 2\nmid k^*\}$. Assume that none of $A\pm 1, 3A\pm 1, A\pm3$ is congruent to $0$ modulo $n$. Then
$$\left|\sum_S \hat{f}(k)\hat{f}(k^*)c_n(Ak+k^*)\right|\le  0.07926\cdot \phi(n).$$
\end{lem}

Remark. The $k^*$ sentence makes sense in view of Lemma~\ref{lem gcd}.

\proof \textit{Case 1.} Both components in every member in $S$ has absolute value $>1$.

By Cauchy-Schwarz, 
\begin{align*}
	\left|\sum_S \hat{f}(k)\hat{f}(k^*)c_n(Ak+k^*)\right|&\le \phi(n)\left(\sum_S |\hat{f}(k)|^2\right)^{1/2}\left(\sum_S |\hat{f}(k^*)|^2\right)^{1/2}\\&\le \phi(n)\cdot \left(\sum_{2\nmid k, k\ne \pm 1} \frac{1}{(k\pi)^2}\right)
	\le \phi(n)\cdot\left(\frac{2}{\pi^2}\cdot(\frac{\pi^2}{8}-1)\right)\\
	&< 0.07926\cdot \phi(n).
\end{align*}

\textit{Case 2.} There exists a member $(k_0,k_0^*)$ in $S$ such that $|k_0|=1$ or $|k_0^*|=1$.

We only prove for  $|k_0|=1$. The case  $|k_0^*|=1$ is similar. So $(1,1^*)\in S$. Thus $0<|1^*|\le H$. Observe that for any $k$ with $0<|k|\le H$, 
$$(Ak+k\cdot 1^*,n)\ge (A+1^*,n)>\sqrt{2H^2n}.$$
Thus, $k\cdot 1^*=k^*$ by the uniqueness of $k^*$.

\textit{Subcase 1.} $(A+1^*,n)=n$.

By assumption, $|1^*|\ne 1,3.$ Therefore, 
\begin{align*}
\left|\sum_S \hat{f}(k)\hat{f}(k^*)c_n(Ak+k^*)\right|&\le \phi(n)\cdot \left(\sum_{2\nmid k} \frac{1}{k\pi\cdot k^*\pi}\right)
\\&\le \phi(n)\cdot\left(\frac{2}{\pi^2}\cdot\frac{\pi^2}{8}\cdot \frac{1}{|1^*|}\right)\\
&\le \phi(n)\cdot\left(\frac{2}{\pi^2}\cdot\frac{\pi^2}{8}\cdot \frac{1}{5}\right)\\
&<  0.07926\cdot \phi(n).
\end{align*}

\textit{Subcase 2.} $(A+1^*,n)<n$.

By the well-known formula for Ramanujan sum, $$|c_n(A+1^*)|\le \phi(n)/\phi\left(\frac{n}{(n,A+1^*)}\right).$$
Since $(n,6)=1$, we have $\phi\left(\frac{n}{(n,A+1^*)}\right)\ge \phi(5)=4.$ Thus $|c_n(A+1^*)|\le \phi(n)/4$. Similarly $|c_n(3A+3^*)|\le \phi(n)/4$ and $|c_n(9A+9^*)|\le \phi(n)/4$. Therefore,
\begin{align*}
&\left|\sum_S \hat{f}(k)\hat{f}(k^*)c_n(Ak+k^*)\right|\le 2|\hat{f}(1)\hat{f}(1^*)c_n(A+1^*)|+2|\hat{f}(3)\hat{f}(3^*)c_n(3A+3^*)|\\&\qquad \qquad \qquad\qquad \qquad \qquad\qquad  + 2|\hat{f}(9)\hat{f}(9^*)c_n(9A+9^*)| + 2\sum_{k\ge 5, k\ne 9, 2\nmid k}|\hat{f}(k)\hat{f}(k^*)|\phi(n)\\
&\qquad\qquad \le 2\cdot \frac{1}{\pi^2}\cdot \frac{\phi(n)}{4}+2\cdot \frac{1}{9\pi^2}\cdot \frac{\phi(n)}{4}+2\cdot \frac{1}{81\pi^2}\cdot \frac{\phi(n)}{4}+2\cdot \frac{1}{\pi^2}\cdot \left(\frac{\pi^2}{8}-1-\frac{1}{9}-\frac{1}{81}\right) \cdot \phi(n)\\
&\qquad \qquad <  0.07926 \cdot \phi(n).
\end{align*}
\qed

\begin{lem}\label{lem 1-3}
	Let $(A,n)=1$. Let $k^*$ and $S$ be the same as in Lemma~\ref{lem *sum}. If any of $3A\pm 1$ or $A\pm3$ is congruent to $0$ modulo $n$. Then
	$$\left|\sum_S \hat{f}(k)\hat{f}(k^*)c_n(Ak+k^*)\right|\le \frac{1}{12}\cdot \phi(n).$$
	(Note that $1/12=0.08333...$.)
\end{lem}

\proof Without loss of generality, we assume $A+3\equiv 0\pmod n.$ Then similar to the proof in Lemma~\ref{lem *sum} we have
\begin{align*}
\left|\sum_S \hat{f}(k)\hat{f}(k^*)c_n(Ak+k^*)\right|
&\le 2\cdot \frac{1}{\pi^2}\cdot \frac{1}{3} \cdot \phi(n)\cdot \left(\sum_{2\nmid k, k\ge 1}\frac{1}{k^2}\right)\\
&= \frac{1}{12} \cdot \phi(n).
\end{align*}
\qed 

\begin{lem}\label{lem s0-s1}
For integers $a, b$ coprime to $n$, define
\begin{align*}
S_0&=\sum_g\chi(g/n)\chi(ag/n)\chi(bg/n)\\
S_1&=\sum_g f(g/n)f(ag/n)f(bg/n).
\end{align*}
We have $$|S_0-S_1|\le \frac{13.02}{H}\cdot \phi(n)+20.02\sqrt{2Hn}+7H$$ as long as $H>1000.$
\end{lem}

\proof  We start with a method appeared in the proof of Erd\H{o}s-Tur\'an-Koksma's Inequality (see~\cite{D-T}, page 21). Writing $x_1=ag/n, x_2=bg/n, x_3=g/n$, we have
\begin{align*}
	|S_0-S_1|&=\left|\sum_g\chi(x_1)\chi(x_2)\chi(x_3)-\sum_g f(x_1)f(x_2)f(x_3)\right|\\
	&\le\sum_g\left(\prod_{j=1}^3(1+|\chi(x_j)-f(x_j)|)-1\right) \\
	& \le \sum_g\left(\prod_{j=1}^3\left(1+\frac{1}{H+1}\sum_{|h_j|\le H}\hat{K}_{H+1}(h_j)C_{h_j}e(h_jx_j)\right)-1\right).
\end{align*}
Next, we open the product, and the above quantity becomes
\begin{align*}
\sum_g\left(\left(1+\frac{1}{H+1}\right)^3-1\right)+\left(1+\frac{1}{H+1}\right)^2\cdot\left(\frac{1}{H+1}\right)\cdot\sum_{0<|h|\le H}\hat{K}_{H+1}(h)C_h\sum_j \sum_g e(hx_j)\\
+\left(1+\frac{1}{H+1}\right)\cdot\left(\frac{1}{H+1}\right)^2\cdot\sum_{h_1}\sum_{h_2}\hat{K}_{H+1}(h_1)C_{h_1}\hat{K}_{H+1}(h_2)C_{h_2}\cdot\\ \sum_g\left(e(h_1x_1+h_2x_2)+e(h_1x_2+h_2x_3)+e(h_1x_3+h_2x_1)\right)\\
+\left(\frac{1}{H+1}\right)^3\cdot \sum_{h_1}\sum_{h_2}\sum_{h_3}\hat{K}_{H+1}(h_1)C_{h_1}\hat{K}_{H+1}(h_2)C_{h_2}\hat{K}_{H+1}(h_3)C_{h_3}\cdot \sum_ge(h_1x_1+h_2x_2+h_3x_3).
\end{align*}

It is straightforward to obtain $\sum_g\left(\left(1+\frac{1}{H+1}\right)^3-1\right)<3.01\phi(n)/H$ as long as $H>1000.$

For the single sum over $h$ we use the well-known bound for Ramanujan sum $|c_n(m)|\le (n,m)$ and the bounds $|\hat{K}|\le 1, |C_h|\le 1$ and the assumption $H>1000$. We see that
\begin{align*}
 \left(1+\frac{1}{H+1}\right)^2\cdot\left(\frac{1}{H+1}\right)\cdot\sum_{0<|h|\le H}\hat{K}_{H+1}(h)C_h\sum_j \sum_g e(hx_j)
\le 7H.
\end{align*}

For the double sum over $h_1,h_2$ we use Lemma~\ref{lem gcd} and get
\begin{align*}
&\left(1+\frac{1}{H+1}\right)\cdot\left(\frac{1}{H+1}\right)^2\cdot\sum_{h_1}\sum_{h_2}\hat{K}_{H+1}(h_1)C_{h_1}\hat{K}_{H+1}(h_2)C_{h_2}\cdot\\ &\qquad \qquad\qquad\qquad\qquad\qquad\sum_g\left(e(h_1x_1+h_2x_2)+e(h_1x_2+h_2x_3)+e(h_1x_3+h_2x_1)\right)\\
&\qquad\le 3\cdot\left(1+\frac{1}{H+1}\right)\cdot\left(\frac{1}{H+1}\right)^2\cdot(2H)\cdot\left(\phi(n)+2H\sqrt{2Hn}\right)\\
&\qquad\le \frac{6.01}{H}\cdot \phi(n)+12.02\sqrt{2Hn}.
\end{align*}

For the triple sum we similarly get
\begin{align*}
\left(\frac{1}{H+1}\right)^3\cdot \sum_{h_1}\sum_{h_2}\sum_{h_3}\hat{K}_{H+1}(h_1)C_{h_1}\hat{K}_{H+1}(h_2)C_{h_2}\hat{K}_{H+1}(h_3)C_{h_3}\cdot \sum_ge(h_1x_1+h_2x_2+h_3x_3)\\
\le \frac{4}{H}\cdot \phi(n)+8\sqrt{2Hn}.
\end{align*}
Collecting the above estimates we are done. \qed

\begin{lem}\label{lem not all 1-3}
Let $(1)(a)(b)(c)$, where $c=2n-1-a-b$, be a minimal zero-sum  sequence modulo $n$, whose index is $2$. Then at least one of the pairs $(1,a),(1,b),(a,b)$ does not satisfy any linear relations $x\pm 3y\equiv 0 \pmod n$ or $3x\pm y\equiv 0\pmod n$ for $(x,y).$
\end{lem}

\proof Suppose for the sake of contradiction that each of the three pairs $(1,a),(1,b),(a,b)$ satisfies a certain linear equation in the lemma. There are $4$ choices of linear equations for each pair, so in total there are $4^3$ possibilities. We claim that none of the pairs could satisfy $3x-y\equiv 0$. Take $(1,a)$ for example. If $3-a\equiv 0\pmod n$, then $a=3$, so $b+c=2n-4$. This implies $b=c=n-2$, and we know such sequence has index $1$. So the claim follows. A similar argument shows that none of the pairs could satisfy $x-3y\equiv 0$. Thus, there are $2^3$ possibilities remaining to be ruled out. These can be done by straightforward verification. To give a flavor, assume $3+a\equiv 1+3b \equiv 3a+b\equiv 0 \pmod n.$ Then we conclude that $28\equiv 0\pmod n$ which implies that $n\le 28$, and we know for such $n$ the index is $1$. \qed

\section{Proof of Theorem~\ref{thm}}
It suffices to prove that
$$S_0=\sum_g\chi(g/n)\chi(ag/n)\chi(bg/n)\ge c\cdot \phi(n)$$ for large $n$. 
Write $S_0=S_1+(S_0-S_1)$, where 
$$S_1=\sum_g f(g/n)f(ag/n)f(bg/n).$$
Using Fourier expansion $f(x)=\sum_{0\le |h|\le H} \hat{f}(h)e(hx)$ where $$ \hat{f}(h)=\hat{\chi}(h)\hat{J}_{H+1}(h),$$ we have
$$S_1=\sum_{h_1}\sum_{h_2}\sum_{h_3}\hat{f}(h_1)\hat{f}(h_2)\hat{f}(h_3)\sum_g e\left(\frac{g}{n}(ah_1+bh_2+h_3)\right).$$
Now $\hat\chi (0)=1/2, \hat\chi (h)=0$ for nonzero even $h$, and that $\hat\chi (h)=1/(\pi i h)$ (in particular, this is purely imaginary) for odd $h$,  and $\hat{J}_{H+1}(h)$ is real. Hence, in the above sum we either have $h_1=h_2=h_3=0$, or exactly one of the three $h$'s equals $0$ (all other terms cancel out). Thus we see that
$$S_1=\phi(n)\cdot(\hat{f}(0))^3+\hat{f}(0)\cdot\left(\sum_{h_2}\sum_{h_3}+\sum_{h_3}\sum_{h_1}+\sum_{h_1}\sum_{h_2}\right),$$
where $$\sum_{h_2}\sum_{h_3}=\sum_{0<|h_2|\le H}\sum_{0<|h_3|\le H} \hat{f}(h_2)\hat{f}(h_3)\sum_g e\left(\frac{g}{n}(bh_2+h_3)\right)$$
and similar  for the other two sums. We now write
$$\sum_{h_2}\sum_{h_3}=S_b^*+T_b$$ where 
$S_b^*$ is the sum in Lemma~\ref{lem *sum} (and Lemma~\ref{lem 1-3}) with $A=b$ in the lemmas, and $T_b$ is a double sum over all the pairs $(h_2,h_3)$ that are not included in $S_b^*$. Using trivial bound for the Ramanujan sums we have
$$|T_b|\le \left(\sum_{0<|h|\le H}\left|\hat f(h)\right|\right)^2\cdot \sqrt{2H^2n}\le \left(\frac{2\log H}{\pi}\right)^2\sqrt{2H^2n}.$$ 
Apply Lemmas~\ref{lem *sum} and~\ref{lem 1-3} we have
$$|S_b^*|\le \frac{1}{12}\cdot \phi(n).$$
We do similar things for $\sum_{h_3}\sum_{h_1}$ and $\sum_{h_1}\sum_{h_2}$. Finally, by assumption  at least one of the three pairs does not satisfy any relation $x\pm 3y\equiv 0 \pmod n$ or $3x\pm y\equiv 0\pmod n$, so at least one of the $S_b^*, S_a^*, S_{ab^{-1}}^*$ is bounded by Lemma~\ref{lem *sum}, that is, $0.07926\cdot\phi(n)$. Thus, collecting all these estimates and keeping in mind that $\hat{f}(0)=1/2$, we obtain that
$$S_1\ge c_0\cdot\phi(n)-\frac{3}{2}\left(\frac{2\log H}{\pi}\right)^2\sqrt{2H^2n},$$
where $$c_0=\frac{1}{8}-\frac{1}{2}\left(0.07926+\frac{2}{12}\right)=0.002....$$

Recall Lemma~\ref{lem s0-s1} says that $|S_0-S_1|\le \frac{13.02}{H}\cdot \phi(n)+20.02\sqrt{2Hn}+7H$  as long as $H>1000.$ It follows that for such $n$
$$S_0\ge \left(c_0-\frac{13.02}{H}\right)\cdot\phi(n)-\frac{3}{2}\left(\frac{2\log H}{\pi}\right)^2\sqrt{2H^2n}-20.02\sqrt{2Hn}-7H.$$
Take $H=13.02\times 10^3,$ say. So $c_0-\frac{13.02}{H}> 0.001=c_1$. Then we see that for some absolute constant $c>0$, $$S_0>c\cdot\phi(n)$$
as long as $c_1\phi(n)-\frac{3}{2}\left(\frac{2\log H}{\pi}\right)^2\sqrt{2H^2n}-20.02\sqrt{2Hn}-7H\ge 0$. A simple computation shows that it suffices to require $\phi(n)/\sqrt{n} > 1.1\times 10^9.$ Using estimates for the $\phi$ function we see that it suffices to take $n>N=10^{20}$. \qed

\section{Proof of Theorem~\ref{thm index}}

We may assume $S=(1)(a)(b)(c)$ where the index of $S$ is $2$ and $(n,ab)=1$ . Also assume $n>N.$ Combine Lemma~\ref{lem not all 1-3} and Theorem~\ref{thm} and we see that there is at least one $g$ (in fact, a positive proportion of $g$'s) makes
$(g)_n/n\in I,(ag)_n/n\in I,$ and $(bg)_n/n\in I$. Thus by Lemma~\ref{lem 3 fall in I} we conclude that the index of $S$ is $1$, a contradiction. \qed

\section*{Acknowledgements}

The author is grateful to Andrew Granville for an insightful conversation which leads to the proof here. This conversation occurred during the conference Probability in Number Theory at CRM in 2018. The author thanks the organizers and the CRM for hosting this wonderful conference. He also thanks Weidong Gao and Alfred Geroldinger for helpful communication, and the anonymous referees for carefully reading the manuscript and making a number of suggestions.

The author is partially supported by a startup fund and a summer research award from William \& Mary, and by a Ralph E. Powe Junior Faculty Enhancement Award from Oak Ridge Associated Universities.

\end{document}